\documentclass[journal]{IEEEtran}
\usepackage{amssymb}
\usepackage{amsthm}

\usepackage{amsmath}
\usepackage{algorithm}
\usepackage{algorithmic}

\usepackage{bm}

\usepackage{mathrsfs}
\usepackage{diagbox}
\usepackage{mathtools}

 \usepackage{color}

\newtheorem{theorem}{Theorem}

\newtheorem{prop}{Proposition}

\newtheorem{dfn}{Definition}

\usepackage{tikz}
\usetikzlibrary{arrows}

\ifCLASSINFOpdf
\else
\fi
\begin{document}
%
\title{A randomized algorithm for the QR decomposition-based approximate SVD}
%
%
%

\author{Xiaohui~Ni,  An-Bao Xu$^{*}$
	\thanks{This work was supported by the National Natural Science Foundation of China under Grant 11801418.}
	\thanks{$^{*}$ Corresponding author.} 
	\thanks{E-mail: 21451025011@stu.wzu.edu.cn (X. Ni), xuanbao@wzu.edu.cn (A.-B. Xu) }}
%
%

\markboth{  }%
{Shell \MakeLowercase{\textit{et al.}}: Bare Demo of IEEEtran.cls for Journals}
%



\maketitle

\begin{abstract}
Matrix decomposition is a very important mathematical tool in numerical linear algebra for data processing. In this paper, we introduce a new randomized matrix decomposition algorithm, which is called randomized approximate SVD based on Qatar Riyal decomposition (RCSVD-QR). Our method utilize random sampling and the OR decomposition to address a serious bottlenck associated with classical SVD. RCSVD-QR gives satisfactory convergence speed as well as accuracy as compared to those state-of-the-art algorithms. In addition, we provides an estimate for the expected approximation error in Frobenius norm. Numerical experiments verify these claims.

\end{abstract}


%
\IEEEpeerreviewmaketitle

\section{Introduction}
%
%
%
%
\IEEEPARstart{M}{atrix} decomposition is a basic tool for data processing in statistics and machine learning \cite{A1}. The most common matrix decomposition methods include singular value decomposition (SVD) \cite{G2}\cite{A2}\cite{Y2}, the rank-revealing QR decomposition (QR) \cite{T3} and the pivoted QLP decomposition (QLP) \cite{G6}. They are widely used in image processing \cite{C4}, signal processing \cite{R5}, information compression \cite{D5} and medicine. In the current era of big data, singular value decomposition is still evident in all types of applications. However, as the size of data matrices generated or collected by computer simulations, experiments, detections and observations continues to increase rapidly, calculating singular value decomposition of large-scale matrices has become a great challenge.

The bottleneck of the traditional method is that its computational complexity is too high, and there are certain limitations on the calculation of large data. In order to improve the speed and accuracy of the algorithm, the pivoted QLP decomposition and the method of computing an approximate SVD based on QR decomposition (CSVD-QR) \cite{L7} are proposed successively. In recent years, low rank approximation randomization algorithm has attracted much attention \cite{G7}\cite{S7}\cite{A7}. Compared with the deterministic algorithm, the randomized low-rank approximation algorithm has the advantages of low complexity, fast running speed and easy implementation. The proposed algorithm provides an idea direction for optimization to obtain more efficient matrix decomposition method.
Based on recent research progress, this paper proposes a new random decomposition method based on CSVD-QR called RCSVD-QR. The calculation of RCSVD-QR requires only matrix-matrix multiplication and several iterative QR decomposition. The introduction of efficient and accurate randomization technology optimizes the algorithm to a certain extent, and is superior to the traditional algorithm in many practical situations.

\subsection{Notations}
In this paper, matrices are represented in uppercase letters.  $   \mathcal{R} (A )$    denotes the numerical range of A, and the i-th largest singular value of $A$ is denoted by $   \sigma _{i}\left (A\right ) $.  The $   \ell_{1}$-norm    and the Frobenius norm of the matrix $   A $    are expressed by $   \left \| A  \right \| _{1} $    and $   \left \| A  \right \| _{F} $.  And $   A ^{\dagger}$    denotes the Moore-Penrose pseudo-inverse of  $A$. In the following, $   \Omega$    represents the random matrix and $   \mathbb{E}$     represents the expected value of the random variable. 
\subsection{Organization}
The rest of this article is structured as follows. 
In section $   \mathrm{II}$,  we describe our proposed new algorithm. In section $   \mathrm{III}$,  we have a theoretical analysis of the algorithm. Finally, we provide numerical experiments to verify the performance of the proposed algorithm in section $   \mathrm{III}$.

\section{RELATED WORK }
Give an $   m\times n$    matrix $A$, we consider a randomized algorithm for computing the largest $r$ ($   r \in( 0,n]$   ) singular values and the corresponding singular vectors of $A$ by QR decompositions directly. More specifically, we can find three matrices $L$,$D$,$R$ such that
\begin{equation}\label{1.1}  A\approx LDR, \end{equation}
where $    L \in R^{m\times r}$    is a column orthogonal matrix, $   D\in R^{r\times r}$    is a diagonal matrix and $    R \in R^{r\times n}$    is a row orthogonal matrix. The diagonal elements of $D$ are equal to the singular value of the SVD-SIM \cite{P8} method. Based on SVD-SIM, the method of calculating the first $r$ singular values and singular vectors by QR decomposition is called CSVD-QR.

The following is a brief description of the CSVD-QR process. Suppose that the results of the $j$-th iteration are $   L_{j}, D_{j}, $    and $   R_{j}$,  with $   L_{0}={\rm eye}(m,r)$,  $   D_{0}={\rm eye}(r,r)$  and $   R_{0}={\rm eye}(r,n)$.  First, $    L_{j}$    is given by QR decomposition of $   AR_{j-1}^{T}$    as
\begin{equation}\label{2}  AR_{j-1}^{T}=L_{j}T \end{equation}
Then, $    R_{j+1}$    is given by the same way
\begin{equation}\label{3}   A^{T}L_{j}=R_{j}D_{j} \end{equation}
Finally, $    D_{j}$    is updated by 
\begin{equation}\label{4}   D_{j}=D_{j}^{T},  \end{equation}
For details, please refer to \cite{L7}.
\section{ALGORITHM MODEL}
In this section, we present a randomized algorithm termed randomized CSVD-QR. Given an matrix $   A\in \mathcal{R} ^{m\times n} $  and an integer $   k \le r<n$,  the basic form of RCSVD-QR is computed as follows: we form a Gauss matrix $   \Omega \in \mathcal{R} ^{n\times r}$  which has independently sampled entries with zero mean and unit variance. First, we compute the following product:
\begin{equation}\label{5}
  Y=A\Omega 
\end{equation}
where $   Y\in \mathcal{R} ^{m\times r} $    is formed by linear combinations of columns of $A$ by the random Gaussian matrix. Then using the QR decomposition algorithm, the factor $Y$ such that : 
\begin{equation}\label{6}
  Y=WC
 \end{equation}
 where $W$ are approximate bases for  $   \mathcal{R} (A)$    and its columns form an orthormal basis for the range of $Y$. We now form a matrix $   B\in \mathcal{R} ^{r\times n} $    as follows: 
 \begin{equation}\label{7}
  B=W^{T}A
\end{equation}
We then compute the CSVD-QR of $B$ with $j$-th iteration: 
 \begin{equation}\label{8}
  B = L_{j} D_{j}^ {T}R_{j}
\end{equation}
Finally, we form the RCSVD-QR-based low-rank approximation of $A$:
 \begin{equation}\label{9}
  A \approx WL_{j} D_{j}^ {T }R_{j}=LDR
\end{equation}

 The pseudo-code for a rank-$k$ RCSVD-QR decomposition algorithm is given in Algorithm 1. 
\begin{algorithm}[htb]
\label{alg:R}  
\caption{Randomized CSVD-QR (RCSVD-QR)}
\begin{algorithmic}[1]
\REQUIRE
Given an matrix $   A\in \mathcal{R} ^{m\times n} $,  a target rank $   k\ge 2$,  an number of inner QR iterations t, an oversampling parameter p, and $   r= k+p$.  Set $   R_{0}={\rm eye}\left ( r, n \right ),  D_{0}= {\rm eye}\left ( r, r \right ) $ 
\ENSURE
 A rank-$k$ approximation\\
\STATE Generate a  Gaussian random matrix  $   \Omega \in \mathcal{R} ^{n\times r}$  \\ 
\STATE Compute a sampling matrix $   Y=A\Omega$  \\ 
\STATE Construct an $   m\times r$    matrix $W$ whose columns form an orthonormal basis for the range of $Y$. (Compute QR factorization $   Y=WC$   ) \\
\STATE Form $   B = W^{T}A$    so that $   P_{W}A=WB$ \\ 
\STATE for $   j=1, 2, \dots, t$    (Compute the CSVD-QR decomposition so that $    B = L_{j} D_{j}^ {T}R_{j}.$   )\\
\STATE Compute QR factorization $   BR_{j-1}^{_T}=L_{j}T$    \\ 
\STATE Compute QR factorization $   B^{T } L_{j}  =R_{j}D_{j}$    and form $   R_{j}=R_{j}^ {T}$   \\
\STATE end\\
\STATE Set $   L=WL_{t}, D=D_{t}^{T}, \ and \ R=R_{t}$   \\
 $   A\approx WL_{t} D_{t}^ {T} R_{t}$   
\end{algorithmic}
\end{algorithm}

\section{ANALYSIS OF RCSVD-QR}
\begin{dfn}\label{definition 2.1}\cite{N9}
	When matrix $P$ is an Hermitian matrix with it satisfies $   P^{2}=P$,  we called it orthogonal projection operator. For a given matrix $M$ with full column rank, this projector can be definitely expressed as follows
	\begin{equation}\label{10}  P_{M}=M(M^{T}M)^{-1}M^{T}. \end{equation}
	where $   M^{T}$    is the transpose matrix of  $M$.   \end{dfn}

Based on (\ref{5}), we have a sample matrix $Y$. Applying Definition 1, the orthogonal projector onto $   \mathcal{R} (Y)$    is 
\begin{equation}\label{11}  P_{Y}=Y(Y^{T}Y)^{-1}Y^{T}. \end{equation}
According to QR decomposition (6), matrix $Y$ can be factorized as $Y=WC$, where  $W$ is a column orthogonal matrix. Then the following formula will be obtained as 
\begin{equation}\label{12}  P_{Y}=Y(Y^{T}Y)^{-1}Y^{T}=WW^{T}. \end{equation}
From the above formula, we can get the following theorem.

\begin{theorem}\label{Theorem 2.2}
	As for RCSVD-QR of  $A\in \mathcal{R} ^{m\times n} $, $A=LDR$, we have the relationship between $\left \| A- LDR \right \| _{F}$ and $\left \| \left( I-P_{Y}\right)A \right \| _{F}$, i.e.,
	\begin{equation}\label{14} 
		\left \| A- LDR \right \| _{F}=\left \| A-WW^{T}A  \right \| _{F}=\left \| \left( I-P_{Y}\right)A \right \| _{F}
	\end{equation}
\end{theorem}
\proof
In order to prove (\ref{14}) , we just have to prove 
\begin{equation}\label{15} 
	LDR = WL_{t} D_{t}^ {T}R_{t}=WB=WW^{T}A
\end{equation}
According to the Algorithm 1 , we know $    LDR = WL_{t} D_{t}^ {T}R_{t}$    and $   WB=WW^{T}A$.  So we only need to prove $   B=L_{t}D_{t}^{T}R_{t}$. 
We will use mathematical induction to prove this equation.

When $   j=1$    , we have 
\begin{equation}\label{16}
	BR_{0}^{T}=L_{1}T
\end{equation}
\begin{equation}\label{17}
	B^{T }L_{1}  =R_{1}D_{1}
\end{equation}
Since $   R_{0}={\rm eye}(r, n)$, $   L^{T}L=I,  RR^{T}=I$,  we have$$   R_{0}^{T}T^{T}L_{1}^{T}L_{1}=T^{T}=R_{1}D_{1}$$ 
$$    BR_{0}^{T}=L_{1}D_{1}^{T}R_{1}^{T}$$ 
Thus $$   B=L_{1}D_{1}^{T}R_{1}^{T}$$
And then $   R_{1}=R_{1}^{T}$, now 
\begin{equation}\label{18}
	B=L_{1}D_{1}^{T}R_{1}
\end{equation}
when $   j=2$  we have
\begin{equation}\label{19}
	BR_{1}^{T}=L_{2}T
\end{equation}
\begin{equation}\label{20}
	B^{T }L_{2}  =R_{2}D_{2}
\end{equation}
Apply (\ref{18}), we have $$   BR_{1}^{T}=L_{1}D_{1}^{T}R_{1}R_{1}^{T}=L_{1}D_{1}^{T}$$ 
So $$   L_{1}D_{1}^{T}=L_{2}T$$
$$   D_{1}=T^{T}L_{2}^{T}L_{1}$$ 
Since $   L_{j}$    is an orthonormal matrix , so$$   B^{T }L_{2}  =R_{1}^{T}D_{1}L_{1}^{T}L_{2}=R_{1}^{T}T^{T}L_{2}^{T}L_{1}L_{1}^{T}L_{2}=R_{1}^{T}T^{T}$$ 
Apply (\ref{20}), 
$$   R_{1}^{T}T^{T}=R_{2}D_{2}$$  
$$   T=D_{2}^{T}R_{2}^{T}R_{1}^{T}$$  
Finally in view of (\ref{19}), $$   BR_{1}^{T}=L_{2}D_{2}^{T}R_{2}^{T}R_{1}^{T}$$ 
Therefore, $$   B=L_{2}D_{2}^{T}R_{2}^{T}$$ 
For $   R_{2}=R_{2}^{T}$  so $$   B=L_{2}D_{2}^{T}R_{2}$$  
Suppose that when $    j=t-1$  the formula 
\begin{equation}\label{21}B = L_{t-1} D_{t-1}^ {T}R_{t-1}
\end{equation} 
is workable. Then we have
\begin{equation}\label{22}
	BR_{t-1}^{T}=L_{t}T
\end{equation}
\begin{equation}\label{23}
	B^{T }L_{t}  =R_{t}D_{t}
\end{equation}
Apply (\ref{21}), we have $$   BR_{t-1}^{T}=L_{t-1}D_{t-1}^{T}R_{t-1}R_{t-1}^{T}=L_{t-1}D_{t-1}^{T}$$ 
So $$   L_{t-1}D_{t-1}^{T}=L_{t}T$$  
$$   D_{t-1}=T^{T}L_{t}^{T}L_{t-1}$$ 
Since $   L_{j}$    is an orthonormal matrix , so$$   B^{T }L_{t}  =R_{t-1}^{T}D_{t-1}L_{t-1}^{T}L_{t}$$ 
$$   =R_{t-1}^{T}T^{T}L_{t}^{T}L_{t-1}L_{t-1}^{T}L_{t}=R_{t-1}^{T}T^{T}$$ 
Apply (\ref{23}), 
$$   R_{t-1}^{T}T^{T}=R_{t}D_{t}$$ 
$$   T=D_{t}^{T}R_{t}^{T}R_{t-1}^{T}$$
Finally in view of (\ref{22}), $$   BR_{t-1}^{T}=L_{t}D_{t}^{T}R_{t}^{T}R_{t-1}^{T}$$
Therefore, $$   B=L_{t}D_{t}^{T}R_{t}^{T}$$ 
For $   R_{t}=R_{t}^{T}$  so $$   B=L_{t}D_{t}^{T}R_{t}$$ 
So equation (\ref{14}) is true.$\Box$

After theorem 1 is proved, we can do some decomposition of $A$ from the equation of $\left \| \left( I-P_{Y}\right)A \right \| _{F}$. 
Since $A$ be given $   m \times n$    matrix with SVD $$   A=U\Sigma V^{T } $$ and fix $   k\ge 0$. Here, $k$ is a fixed number.  Roughly,  the subspace spanned by the first $k$ left singular vectors is the target that the proto-algorithm attempts to approximate. For better analysis, the SVD decomposition above is divided as follows \cite{N9}:
\begin{equation}\label{24} 
	A=U \begin{bmatrix}
		\Sigma_{1}  &  \\
		&\Sigma _{2} 
	\end{bmatrix}
	\begin{bmatrix}
		V_{1}^{T }  \\V_{2}^{T }
	\end{bmatrix} 
\end{equation}
The matrices $   \Sigma_{1}\in \mathcal{R} ^{k\times k}$    and $   \Sigma_{2}\in \mathcal{R} ^{(n-k) \times (n-k)}$    are square. In the analysis, we will see that the left unitary factor $U$ does not play a significant role. So we mainly analyze the right unitary factor $V$.

Let $   \Omega$    be an $   n \times r$    test matrix, where $r$ denotes the number of samples. We assume only that $   r \ge k$. In the coordinate system determined by the right unitary factor of $A$, the test matrix is decomposed as follows:
\begin{equation}\label{25} 
	\Omega_{1}=V_{1}^{T }\Omega \quad and \quad   \Omega_{2}=V_{2}^{T }\Omega.
\end{equation}
This equation along with the following Proposition and Theorem, both are the powerful tools used to prove Theorem 3.
\begin{theorem}\label{Theorem 2.3}\cite{C10}(deterministic error bound) Assuming that $   \Omega_{1}$  has full row rank, the approximation error satisfies
	\begin{equation}\label{26} 
		\left \| \left( I-P_{Y}\right)A \right \| _{F}^{2} \le\left \| \Sigma _{2}  \right \|_{F}^{2} + \left \| \Sigma _{2} \Omega _{2} \Omega _{1}^{\dagger }  \right \|_{F}^{2}  
	\end{equation}
\end{theorem}

\begin{prop}\label{Proposition 1 }\cite{Y11}\cite{Y12} (scaling of the expected norm of a Gaussian matrix)
	Fix matrices $S$, $T$, and construct a standard Gaussian matrix $H$. Then 
	\begin{equation}\label{27} 
		\left (E  \left \|  SHT \right \|_{F}^{2}   \right )^{1/2}= \left \| S  \right \|_{F}  \left \| T  \right \|_{F} 
	\end{equation}
\end{prop}
\begin{prop}\label{Proposition 2}\cite{R13}\cite{Z14}(expected norm of a pseudoinverted Gaussian matrix). Construct a $   k\times (k+p)$  standard Gaussian matrix $H$ with $   p\ge 2$  and $   k\ge 2$  Then 
	\begin{equation}\label{29} 
		\left ( \mathbb{E}  \left \| H ^{\dagger }  \right \|^{1/2}  \right ) =\sqrt{\frac{k}{p-1}} 
	\end{equation}
\end{prop}
\begin{theorem}\label{Theorem 2.4} With the notation of  Algorithm 1, choose an oversampling parameter $   p\ge 2$ and a target rank $   k\ge 2$ .Here limit $  k+p=\min(m,n)$. Construct the orthonormal matrix $W$ according to the sample matrix $Y$ and perform the CSVD-QR decomposition on the reduced matrix $   B=W^{T}A$  such that $    B = L_{t} D_{t}^ {T}R_{t}.$   Then the expected approximation error                                                                                                                                                                                                                                                 
	\begin{equation}\label{31} 
		\mathbb{E}\left [\left \| A- LDR \right \| _{F}\right] \le \left ( 1+\frac{k}{p-1}  \right ) ^{1/2}\left (\sum^{n}_{j=k+1}\sigma _{j}^{2}(A)   \right )^{1/2} 
	\end{equation}
	where $   \mathbb{E}\left [ \cdot  \right ]$    denotes expectation, $   L=WL_{t}, D=D_{t}^{T}$,  and $ R=R_{t}$ 
\end{theorem}

\noindent \proof According to the Theorem 1, we know
$$   
\left \| A- LDR \right \| _{F}=\left \| A-WW^{T}A  \right \| _{F}=\left \| \left( I-P_{Y}\right)A \right \| _{F}
$$   
In view of above formula, it suffices to prove that
\begin{equation}\label{32} 
	\mathbb{E}\left [\left \| \left( I-P_{Y}\right)A \right \| _{F}\right] \le \left ( 1+\frac{k}{p-1}  \right ) ^{1/2}\left (\sum^{n}_{j=k+1}\sigma _{j}^{2}(A)   \right )^{1/2} 
\end{equation}
Theorem 2 imply that the right unitary factor $V$ of $A$ is divided into $   V_{1}$    and $   V_{2}$  respectively, $k$ and $   n-k$    columns.  The Gaussian distribution is rotationally invariant, so $   V^{T}\Omega$   is also a standard Gaussian matrix. Since that $    \Omega _{1}$    and $    \Omega _{2}$    are nonoverlapping submatrices of $   V^{T}\Omega$,  so $    \Omega _{1}$    and $    \Omega _{2}$    are not only standard Guassian but also stochastically independent. Moreover, the rows of a (fat) Gaussian matrix are almost certainly in general position, so the $   k \times (k+p)$    matrix $    \Omega _{1}$     has row full rank probability one.

Combine H$   \ddot{o}$lder’s inequality and Theorem 2, we have
\begin{align*} 
	\mathbb{E}\left \| \left( I-P_{Y}\right)A \right \| _{F}  \le {\mathbb{E}\left [\left \| \left( I-P_{Y}\right)A \right \| _{F}^{2}\right]}^{1/2}\\
	\le (\left \| \Sigma _{2}  \right \|_{F}^{2} +\mathbb{E} \left \| \Sigma _{2} \Omega _{2} \Omega _{1}^{\dagger }  \right \|_{F}^{2}  )^{1/2} 
\end{align*}
We compute this expectation by conditioning on the value of $    \Omega _{1}$    and applying Proposition 1  to the scaled Gaussian matrix $    \Omega _{2}$    . Thus, 
\begin{align*} 
	\mathbb{E} \left \| \Sigma _{2} \Omega _{2} \Omega _{1}^{\dagger }  \right \|_{F}^{2} 
	&=\mathbb{E}\left (  \mathbb{E} \left [   \left \| \Sigma _{2} \Omega _{2} \Omega _{1}^{\dagger }  \right \|_{F}^{2}\mid \Omega _{1} \right ]\right ) \\
	&=\mathbb{E} \left ( \left \| \Sigma _{2}  \right \| _{F}^{2}  \left \| \Omega _{1}^{\dagger }  \right \| _{F}^{2} \right ) \\
	&=\left \| \Sigma _{2}  \right \| _{F}^{2}\cdot  \mathbb{E} \left \| \Omega _{1}^{\dagger }  \right \| _{F}^{2}\\
	&=\frac{k}{p-1}\cdot \left \| \Sigma _{2}  \right \| _{F}^{2}
\end{align*}
where the last expectation follows from relation (\ref{29}) of Proposition 2. Thus
\begin{align*} 
	\mathbb{E}\left [\left \| A- LDR \right \| _{F}\right] =\mathbb{E}\left \| \left( I-P_{Y}\right)A \right \| _{F}  \\
	\le (1+\frac{k}{p-1})^{1/2} \left \| \Sigma _{2}  \right \| _{F}
\end{align*}
where  $    \left \| \Sigma _{2}  \right \| _{F}^{2}= \sum^{n}_{j=k+1}\sigma _{j}^{2}(A) $.$\Box$
\\
\section{Numerical Experiments}

To investigate the performance of RCSVD-QR , we tested the following numerical experiments. All tests are performed on the MATLAB R2021b with Intel(R) Core(TM) i5-12500 3.00GHz and 16GB of RAM. 

We compare the performance of RCSVD-QR with the following five types of algorithms, including the SVD, RSVD \cite{N9}, CSVD-QR and Rand-QLP \cite{N15}.

\subsection{Precision Comparison}
\

In this section, we test the convergence of the RCSVD-QR algorithm on synthetic data. The synthetic data in the experiments is generated according to the following formula: 
 \begin{equation}\label{33} 
M=M_{1}^{m\times r_{1}}M_{2}^{ r_{1}\times n}
 \end{equation}
  \begin{equation}\label{34} 
M_{1}^{m\times r_{1}}={\rm randn}(m, r_{1})
 \end{equation}
  \begin{equation}\label{35} 
M_{2}^{ r_{1}\times n}={\rm randn}(r_{1}, n)
 \end{equation}
where $    r_{1}\in \left [1,n\right] $    is the rank of M. 
At first, we calculate the relative error of each iteration of the RCSVD-QR by the formula :
 \begin{equation}\label{36} 
e_{k} =\frac{ {\textstyle \sum_{i=1}^{t}\left |\left \| d_{i}  \right \| _{1} -\sigma _{i}\left ( M \right )   \right | } }{ {\textstyle \sum_{i=1}^{t}}\sigma _{i}\left ( M \right ) }  (t={\rm min}(r, r_{1}))
 \end{equation}
where $   \sigma _{i}\left ( M \right ) $   is the ith singular value of $M$ and $   d_{i}$    is diagonal elements of $   D_{k}$    in the kth iteration of the RCSVD-QR. Here, we set $   m=n=r=300$  and $r_{1}=250$.  

Figure 1 shows that under this set of parameters, the relative error of RCSVD-QR decreases rapidly in the first 5 iterations, and then gradually converges. In order to measure the optimal performance of this algorithm, we set the number of iterations as 5 in the following numerical tests.

\begin{figure}[htbp]
\centering
\includegraphics[width=2.8in]{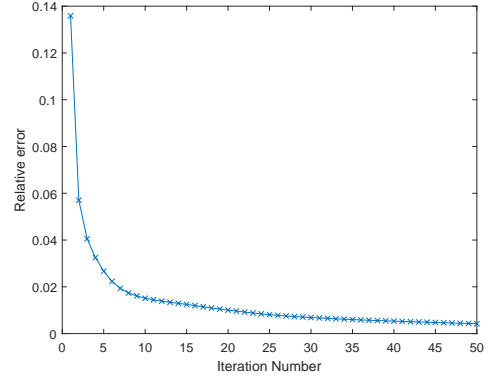}
\caption{Relative error convergence }
\end{figure}

Then we compare the performance of the RCSVD-QR on two kinds of synthetic matrices to verify that the algorithm provides high precision singular value and low rank approximation. 
\begin{itemize}
\item[1)]
A Noisy Low-Rank Matrix \cite{J15}: The first matrix takes the following form as $   A=U\Sigma V^{T}+0.1\sigma_{k}G_{n}$,  where $   U, V\in R^{n\times n}$    are random orthonormal matrices and $   G_{n}$    is a normalized Gaussian matrix. The diagonal elements of $   \Sigma$    are singular values $   \sigma _{i} s$    decreasing linearly from 1 to $   10^{-9}$  and when $   j=k+1,\dots ,n$,  we set $   \sigma _{j} =0$.  Here, we set $   n=10^{3}$    and rank $   k=20$ 
\item[2)]
A Matrix with Rapidly Decaying Singular Values\cite{M16}: The matrix is formed as $   A=U\Sigma V^{T}$  where $U$,$V$ are same as above and $   \Sigma = {\rm diag} \left ( 1^{-1},2^{-1},\dots,n^{-1}   \right ) $.  Here, we set $   n=10^{3}$    and $   k=10$ 
\end{itemize}

The results of singular values calculated by each algorithm for the above two kinds of matrices can be clearly seen in Figure 1 and Figure 2.

\begin{figure}
  \begin{center}
  \includegraphics[width=2.8in]{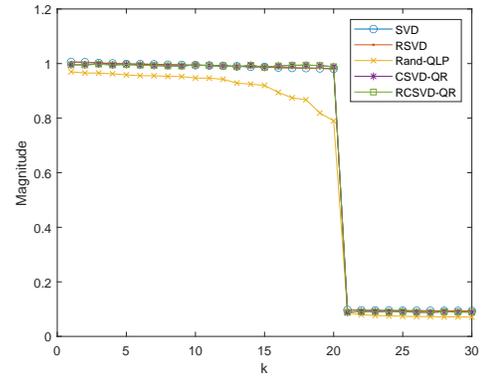}
  \caption{Comparison of singular values for the noisy low-rank matrix. }
  \end{center}
\end{figure}
\begin{figure}
  \begin{center}
  \includegraphics[width=2.8in]{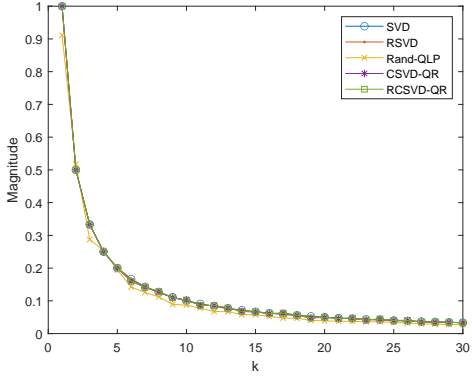}
  \caption{Comparison of singular values for the matrix with Rapidly Decaying Singular Values.}
  \end{center}
\end{figure} 

For the first matrix, the approximation of most of the singular values is a good estimate of the true singular values. In this  is superior to Rand-QLP , and has similar performance to CSVD-QR. Compared with the optimal SVD, the approximate accuracy of RCSVD-QR is higher.

For the second matrix with rapidly decaying singular values, Rand-QLP performed the worst and the other algorithms performed similarly. In contrast, RSVD, CSVD-QR and RCSVD-QR approximate the singular values of the matrix with very high precision.

\subsection{Runtime Comparison}
\

We compare the speed of RCSVD-QR with other 4 kinds  algorithms. Here, we set $   m=n=k(\times1000)$,  $   r_{1}=1000$    and $k$=1 to 10. These algorithms are applied to the synthetic matrix $M$ as follows: 
 Figures 4 and 5 provide CPU time for different algorithms and speedups for RCSVD-QR.

\begin{figure}
  \begin{center}
  \includegraphics[width=2.8in]{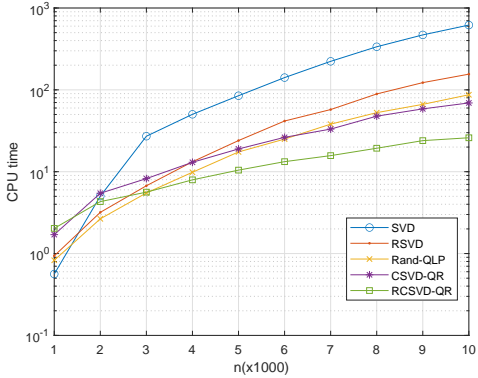}
  \caption{Computational times for different algorithms on random matrices of different sizes.}
  \end{center}
\end{figure}
\begin{figure}
  \begin{center}
  \includegraphics[width=2.8in]{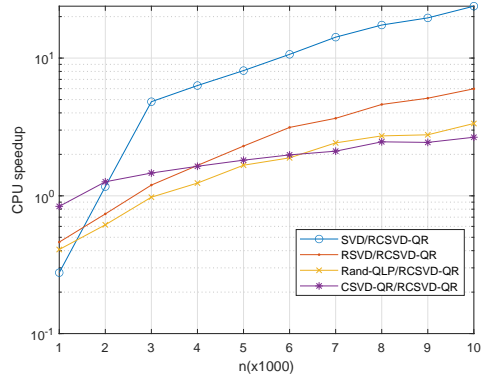}
  \caption{Speedups provided by RCSVD-QR on computation time for random matrices of different scale sizes}
  \end{center}
\end{figure} 

When the random matrix size is larger than $   3000\times 3000$  the speed of RCSVD-QR is the fastest among other algorithms. It can be seen that the speed advantage of this algorithm is more obvious in large-scale matrix. Compare with SVD, RSVD, Rand-QLP and CSVD-QR, the acceleration of RCSVD-QR is up to $    23\times, 6\times, 3\times, 2\times$.  
\subsection{Image Compression}
\

Image compression is a classic application of singular value decomposition. We transform the picture into a matrix, keeping only the first $k$ singular values to achieve the effect of compression size.
\begin{figure}[!t]
\flushleft
\includegraphics[width=3.3in]{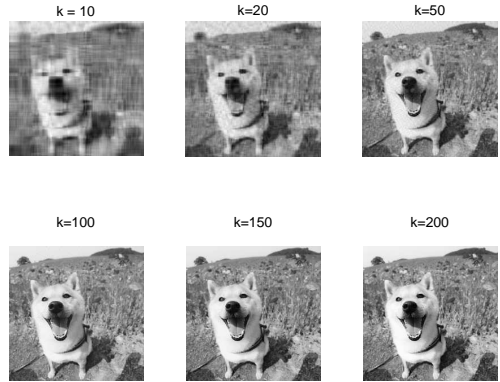}
\caption{Image compression effect with the first $k$ singular values}
\end{figure}
The image compression effect can be seen from Figure 6, when $k$=10, 20, 50, 100, 150 and 200 with RCSVD-QR algorithm. At this time, the number of iterations taken is 5. It can be seen that the image compression effect of this algorithm can be comparable to that of SVD when only 5 iterations are taken, which reflects the high accuracy of this algorithm. When $k$=100, the picture is basically the same as the original picture. Other pictures also achieve better clarity when the first 100 singular values are obtained, so we fixed $k$=100 below. We randomly selected 6 pictures of different sizes for image compression processing with the size of $   512\times 512, 1600\times1600, 1920\times2880, 3456\times5184, 6713\times5168, 5461\times8192$.

Table 1 shows the CPU time compressed by different algorithms for each of the six images. In the experiment of real data, we can see the speed advantage of this algorithm. And the larger the matrix size, the more obvious the advantage. 
\begin{table}[]
\centering
\caption{CPU time of different image compression}
\renewcommand\arraystretch{1.5}
\begin{tabular}{cccccc}
\hline
       & SVD      & RSVD   & Rand  & CSVD & RCSVD  \\ 
       &&&-QR&-QLP&-QR\\\hline
Lena   & 0.1294   & 0.0522 & 0.0494 & 0.0431  & \textbf{0.0416}    \\ 
Dog    & 2.0653   & 1.1297 & 0.8756 & 0.3063  & \textbf{0.1597}    \\ 
Flower & 6.5969   & 3.5172 & 1.6938 & 1.0938  & \textbf{0.2906}    \\ 
Tea    & 61.6125  & 19.3734  & 7.0531  & 2.7641  & \textbf{0.5938} \\ 
Kid    & 181.5391 & 52.7750 & 44.9609 & 7.0359  & \textbf{1.5969}  \\ 
Tree   & 224.1172 & 89.3047 & 25.8016& 7.8594  & \textbf{1.5453}   \\ \hline
\end{tabular}
\end{table}

\section{conclusion}
In this paper, we proposed RCSVD-QR, a randomized rank-revealing algorithm that approximates SVD of a matrix. Then the deviation bound of the algorithm is given on the probability. Numerical experimental results show that this method is not only much fast than the stat-of-the-art algorithm, but it also has good precision.


%

\ifCLASSOPTIONcaptionsoff
  \newpage
\fi

\end{document}